\documentclass[12pt,a4paper]{article}
\usepackage{amssymb}

\usepackage{graphicx,amssymb,amsfonts,epsfig,amsthm,a4,amsmath,url}
\usepackage[latin1]{inputenc}

\psfigdriver{dvips}

\newtheorem{thm}{Theorem}
\newtheorem{cor}[thm]{Corollary}
\newtheorem{lem}[thm]{Lemma}

\newtheorem{prop}[thm]{Proposition}
\newtheorem{clai}[thm]{Claim}
\theoremstyle{definition}
\newtheorem{defn}[thm]{Definition}
\newtheorem{que}[thm]{Question}
\theoremstyle{remark}
\newtheorem{rem}[thm]{Remark}

\newtheorem{exes}[thm]{Examples}

\numberwithin{equation}{section}

\newcommand{\Z}{\mathbf{Z}}

\newcommand{\N}{\mathbf{N}}

\newcommand{\R}{\mathbf{R}}

\newcommand{\tpr}{\begin{tiny}\noindent Proof:}

\newcommand{\bpr}{\noindent \textbf{Proof}: ~}
\newcommand{\epr}{~$\blacksquare$}

\title{Volume of spheres in doubling metric measured spaces and in groups of polynomial growth}%
\author{Romain Tessera }
\date{\today}

\begin{document}

\baselineskip=16pt

\maketitle

 \maketitle
\begin{abstract}
Let $G$ be a compactly generated locally compact group and let $U$
be a compact generating set. We prove that if $G$ has polynomial
growth, then $(U^n)_{n\in \N}$ is a F\o lner sequence and we give
a polynomial estimate of the rate of decay of
$$\frac{\mu(U^{n+1}\smallsetminus U^n)}{\mu(U^n)}.$$
Our proof uses only two ingredients: the doubling property and a
weak geodesic property that we call Property (M). As a matter of
fact, the result remains true in a wide class of doubling metric
measured spaces including manifolds and graphs. As an application,
we obtain a $L^p$-pointwise ergodic theorem ($1\leq p<\infty$) for
the balls averages, which holds for any cglc group $G$ of polynomial
growth.
\end{abstract}


\section{Introduction}
Let $G$ be a compactly generated, locally compact (cglc) group
endowed with a left Haar measure $\mu$. Recall that a sequence
$(A_n)_{n\in \N}$ of subsets of a locally compact group $G$ is
said to be F\o lner if for any compact set $K$,
$$\mu(K.A_n\vartriangle A_n)=o(\mu(A_n)).$$
Let $U$ be a compact generating set of $G$ (we mean by this that
$\cup_{n\in \N} U^n=G$), non necessarily symmetric. If $\mu(U^n)$
grows exponentially, it is easy to see that the sequence
$(U^n)_{n\in \N}$ {\it cannot} be F\o lner. On the other hand, if
$\mu(U^n)$ grows subexponentially, then there exists trivially a
sequence $(n_i)_{i\in \N}$ of integers such that $(U^{n_i})_{i\in
\N}$ is F\o lner. But it is not clear whether the whole sequence
$(U^n)_{n\in \N}$ is F\o lner. This was first conjectured for
amenable groups by Greenleaf in 1969 (\cite{Green}, p 69), who also
proved it with Emerson \cite{EmGr} in the Abelian case, correcting a
former proof of Kawada \cite{Kawa} (see also Proposition
\ref{abelian}). The conjecture is actually not true for all finitely
generated amenable groups since there exist amenable groups with
exponential growth (for instance, all solvable groups which are not
virtually nilpotent). Nevertheless, the conjecture is still open for
groups with subexponential growth. In 1983, Pansu \cite{Pan} proved
it for nilpotent finitely generated groups. In \cite{Br}, Breuillard
recently generalized the theorem of Pansu, which now holds for
general cglc groups of polynomial growth. In fact, They prove that
$\mu(U^n)\sim Cn^d$, for a constant $C=C(U)>0$, which clearly
implies that $(U^n)_{n\in \N}$ is F\o lner. In this article, we
prove the conjecture for all compactly generated groups with
polynomial growth. More precisely, we prove the following theorem:
there exist $\delta>0$ and a constant $C<\infty$, such that
$$\mu\left(U^{n+1}\smallsetminus U^n\right)\leq Cn^{-\delta}\mu(U^n).$$
Interestingly, our proof works in a much more general setting.
Recall that a metric measure space $(X,d,\mu)$ satisfies the
doubling condition (or ``is doubling") if there exists a constant
$C\geq 1$ such that
$$\forall r>0, \forall x\in X,\quad \mu(B(x,2r))\leq C\mu(B(x,r))$$
where $B(x,r)=\{y\in X, d(x,y)\leq r\}$. Let $S(x,r)$ denote the
``$1$-sphere" of center $x$ and radius $r$, i.e.
$S(x,r)=B(x,r+1)\smallsetminus B(x,r)$. Actually, we prove a similar
result for doubling metric measured spaces satisfying a weak
geodesic property we will call Property (M) (see §~\ref{sec}). In
this setting, the result becomes: there exist $\delta>0$ and a
constant $C<\infty$, such that
$$\forall x\in X, \forall r>0, \quad \mu(S(x,r))\leq Cr^{-\delta}\mu(B(x,r)).$$
In particular, the conclusion of this theorem holds for metric
graphs and Riemannian manifolds satisfying the doubling condition.

In the case of metric measured spaces, our result is somewhat
optimal, since in \cite[Theorem~4.9]{tess'}, we build a graph $X$,
quasi-isometric to $\Z^2$, such that there exist $0<a<1$, an
increasing sequence of integers $(n_i)_{i\in\N}$ and $x\in X$ such
that
$$|S(x,n_i)|\geq c|B(x,n_i)|/n_i^a \quad \forall i\in\N.$$
Note that easier counter examples can be obtained with trees with linear
growth (see Remark~\ref{counterexampleRemark}).
Moreover, we will see that our assumptions on $X$, that is, Doubling
Property and Property (M) (see Definition \ref{M} below) are also
optimal in some sense. 

An interesting and historical motivation (see for instance
\cite{Green}) for finding F\o lner sequences in groups comes from
ergodic theory. As a consequence of our result, we obtain a
$L^p$-pointwise ergodic theorem ($1\leq p<\infty$) for the balls
averages, which holds for any cglc group $G$ of polynomial growth
(see theorem \ref{thergo}). We refer to a recent survey of A. Nevo
\cite{Amos} for more details and complete proofs.

\bigskip
\begin{center}
{\bf Note}
\end{center}

After this paper was accepted for publication, Assaf  Naor and  Emmanuel
Breuillard pointed to me a reference containing our main result in the
context of metric measure spaces, namely Theorem~\ref{th1} (see
\cite{CM}[Lemma 3.3]).  

\section{Main results}

\subsection{Property (M)}

\begin{defn}\label{M}
We say that a metric space $(X,d)$ has Property (M) if there exists
$C<\infty$ such that the Hausdorff distance between any pair of
balls with same center and any radii between $r$ and $r+1$ is less
than $C$. In other words, $\forall x\in X$, $\forall r>0$ and
$\forall y\in B(x,r+1)$, we have $d(y,B(x,r))\leq C.$
\end{defn}

\begin{prop}\label{Mprop}
Let $(X,d)$ be a metric space. The following properties are
equivalent
\begin{enumerate}
\item $X$ has Property (M). \item $X$ has ``monotone\footnote{This is why we call this property (M).} geodesics"
, i.e. there exists $C<\infty$ such that, for all $x,y\in X$,
$d(x,y)\geq 1$, there exists a finite chain $x_0=x,
x_1,\ldots,x_m=y$ such that for $0\leq i<m$,
$$d(x_i,x_{i+1})\leq C;$$
and
$$d(x_{i},x)\leq d(x_{i+1},x)-1.$$
\item There exists a constants $C<\infty$ such that for all $r>0$, $s\geq 1$ and $y\in B(x,r+s)$,
$$d(y,B(x,r))\leq Cs.$$
\end{enumerate}
\end{prop}
\noindent{\bf Proof of $1\Rightarrow 2$.} Let $x,y\in X$ be such
that $d(x,y)\geq 1$. Let us construct the sequence
$x=x_0,x_1,\ldots,x_m=y$ inductively. First, by Property (M) and
since $d(x,y)\geq 1$, there exists $x_1\in B(y,d(x,y)-1)$ such that
$d(x,x_1)\leq C$. Now, assume that we have constructed a sequence
$x=x_0,x_1,\ldots,x_k$ such that for $0\leq i<k$,
$$d(x_i,x_{i+1})\leq C;$$
and
$$d(x_{i},x)\leq d(x_{i+1},x)-1.$$
If $d(x_k,y)< 1$, then up to replace $C$ by $C+1$, and $x_k$ by $y$,
the sequence $x_0,\ldots,x_k$ is a monotone geodesic between $x$ and
$y$. Otherwise, there exists $x_{k+1}\in B(y,d(y,x_k)-1)$ such that
$d(x_k,x_{k+1})\leq C.$ Clearly this process has to stop after at
most $[d(x,y)]$ steps, so we are done.

\noindent{\bf Proof of $2\Rightarrow 3$.} Let $x_0=x,
x_1,\ldots,x_m=y$ be a monotone geodesic from $x$ to $y$. There
exists an integer $k\leq s+1$ such that $x_{m-k}\in B(x,r)$. Hence
$$d(y,B(x,r))\leq d(y,x_{m-k})\leq C(s+1)\leq 2Cs$$
which proves the implication.

\noindent{\bf Proof of $3\Rightarrow 1$.} Just take $s=1$.\epr

\

\noindent{\bf Invariance under Hausdorff equivalence.} Recall (see
\cite{Grom}, p.2) that two metric spaces $X$ and $Y$ are said
Hausdorff equivalent
$$X\sim_{Hau}Y$$
if there exists a (larger) metric space $Z$ such that $X$ and $Y$
are contained in $Z$ and such that
$$\sup_{x\in X}d(x,Y)<\infty$$
and
$$\sup_{y\in Y}d(y,X)<\infty.$$

It is easy to see that Property (M) is invariant under Hausdorff
equivalence. But on the other hand, Property (M) is unstable under
quasi-isometry. To construct a counterexample, one can
quasi-isometrically embed $\R_+$ into $\R^2$ such that the image,
equipped with the induced metric does not have Property (M):
consider a stairway-like curve starting from $0$ and containing for
every $k\in \N$ a half-circle of center $0$ and radius $2^k$. So (M)
is strictly stronger than the quasi-geodesic property (\cite{Grom},
p.7), which is invariant under quasi-isometry: $X$ is quasi-geodesic
if there exist two constants $d>0$ and $\lambda>0$ such that for all
$(x,y)\in X^2$ there is a finite chain of points of $X$
$$x=x_0,\ldots,x_m=y,$$
such that
$$d(x_{i-1},x_{i})\leq d,\quad i=1\ldots m,$$
and
$$\sum_{i=1}^m d(x_{i-1},x_i)\leq \lambda d(x,y).$$

\begin{exes}
Recall that a metric space $(X,d)$ is called $b$-geodesic if for any
$x,y\in X$, there exists a finite chain $x=x_0,\ldots,x_m=y$ such
that
$$d(x_i,x_i+1)\leq b \quad \forall i=0,\ldots,m-1$$
and
$$d(x,y)=\sum_{i=1}^m d(x_{i-1},x_i).$$
Note that a $b$-geodesic space is $b'$-geodesic for any $b'\geq b$.
For example, Riemannian manifolds and more generally geodesic spaces
are $b$-geodesic for any $b>0.$ Clearly, $1$-geodesic spaces satisfy
Property (M). Other examples of $1$-geodesic spaces are graphs.
Namely, to any connected simplicial graph, we associate a metric
space, whose elements are the vertices of the graph, the metric
being the usual shortest path distance between two points. We simply
call such a metric space a graph. By definition, graphs are
$1$-geodesic, so in particular they satisfy Property (M). Finally, a
discretisation (i.e. a discrete net) of a Riemannian manifold $M$
has Property (M) for the induced distance (this is a consequence of
the stability under Hausdorff equivalence).
\end{exes}

\subsection{The main theorem}

Let $X=(X,d,\mu)$ be a metric measured space. By metric measure
space, we mean that $\mu$ is a Borel measure, supported on the
metric space $(X,d)$ satisfying $\mu(B(x,r)<\infty$ for all $x\in X$
and $r>0$. Recall that $X$ is said to be doubling if there exists
$C\geq 1$ such that
$$\mu(B(x,2r))\leq C\mu(B(x,r))\quad \forall x\in X, \forall r>0.$$

Our main result says that in a doubling space with Property (M),
balls are F\o lner sets.

\begin{thm}\label{th1}
Let $X=(X,\mu,d)$ be a doubling metric measured space with Property
(M). Then, there exists $\delta>0$ and a constant $C<\infty$ such
that, $\forall x\in X$ and $\forall n\in \N$
$$\mu\left(S(x,n)\right)\leq Cn^{-\delta}\mu(B(x,n)).$$
In particular, the ratio $\mu(B(x,n+1)\smallsetminus
B(x,n))/\mu(B(x,n))$ tends to $0$ uniformly in $x$ when $n$ goes
to infinity.
\end{thm}

\subsection{Optimality of the assumptions}

\noindent{\bf The doubling property.} First, note that the doubling
assumption cannot be replaced by polynomial growth. Indeed, for
every integer $n$, consider the following finite rooted tree: first,
take the standard ternary rooted tree of depth $n$. Then stretch it
as follows: replace each edge connecting a $k-1$'th generation
vertex to a $k$'th generation vertex by a (graph) interval of length
$2^{n-k}$. We obtain a rooted tree $G_n$ of dept
$2^{n}=\sum_{k=0}^{n-1}2^k$. Then consider the graph $G_n'$ obtained
by taking two copies of $G_n$ and identifying the vertices of last
generation of the first copy with those of the second copy. Write
$r_n$ and $r'_n$ for the two vertices of $G_n'$ corresponding to the
respective roots of the two copies of $G_n$. Finally, glue
``linearly" the $G_n'$ together identifying $r_n'$ with $r_{n+1}$,
for all $n$: it defines a infinite connected graph $X$.

Let us prove that this graph has polynomial growth. It is enough to
look at radii of the form $r=2^k$ where $k\in \N.$ On the other hand,
among the balls of radius $2^k$, those which are centered in points
of $n$'th generation of a $G_n$ for $n$ large enough are of maximal
volume. Let us take such an $x$.  Note that $B(x,2^k)$ is isometric
to $G'_{k-1}$ glued with to segments of length $2^{k-1}$ at its
extremities $r_{2^{k-1}}$ and $r'_{2^{k-1}}.$ Hence,
\begin{eqnarray*}
|B(x,2^k)|& = & 2\left(2^{k-1}+\sum_{j=0}^{k-2}3^j2^{k-j}\right)\\
& \leq & 2^k+2^{k+2}(3/2)^{k-1}\\
& \leq & 83^k \\
& = & 8r^{\log 3/\log 2}.
\end{eqnarray*}

On the other hand, the sphere
$S(r_{2^k},2^k)=B(r_{2^k},2^k)\smallsetminus B(r_{2^k},2^{k-1})$ has
a volume larger than $3^k$, so that
$$\frac{|S(r_{2^k},2^k)|}{|B(r_{2^k},2^k)|}\geq 1/8.$$

\begin{rem}\label{counterexampleRemark}
We can generalize the above construction taking trees of valence $b\geq 2$
and replacing each edge connecting a $k-1$'th generation
vertex to a $k$'th generation vertex by a (graph) interval of length
$a^{n-k}$, for some $a\in \N$. If $a<b$, then we still obtain graph with
polynomial growth, but not doubling which contradicts Theorem~\ref{th1}.
But if $a>b$, then it is easy to see that we obtain a graph with linear
growth, hence doubling. Moreover, there exists a sequence of vertices
$x_n$ and $c>0$ such that
$$|S(x_n,a^n)|\geq c b^n.$$
Hence, in this case the $\delta$ of Theorem~\ref{th1} is less than $1-\log
b/\log a$.
\end{rem}

\noindent{\bf The property (M).} Another interesting point is the
fact that Property (M) cannot be replaced by any quasi-isometry
invariant property like quasi-geodesic property.

Indeed, one can very easily build a counterexample, embedding
quasi-isometrically $\R_+\times [0,1]$ into $\R^2$. First, consider
a stairway-like curve starting from $0$ and containing for every
$k\in \N$ a half-circle of radius $2^k$ (c.f. the counterexample to
Property (M)), and then look at the closed $1$-neighborhood of this
curve. Denote by $X$ the corresponding closed subset of $\R^2$,
equipped with the Lesbegues measure and the induced distance. As $X$
is quasi-isometric to $\R_+$, the volume of balls grows linearly.
But observe that the volume of $B_X(0,2^k+1)\smallsetminus
B_X(0,2^k)$ is larger than $\pi 2^k$.

In particular, balls being F\o lner sets is not invariant under
quasi-isometry.

\subsection{An interesting particular case: locally compact groups with polynomial
growth}\label{group}

Let $(G,\mu)$ be a cglc group endowed with a Haar measure $\mu$.
Let $U$ be a compact generating set of $G$. Define a left
invariant distance $d$ on $G$ by:
$$\forall x,y\in G, \quad d(x,y)=\inf\{n\in \N, x^{-1}y\in U^n\}.$$ Note that unless $U$ is symmetric (i.e.
$U^{-1}=U$), $d$ is not really a distance since we do not have:
$d(x,y)=d(y,x)$. Nevertheless,  $d$ is ``weakly" symmetric, i.e.
there exists a constant $C<\infty$ such that
$$\forall x,y\in X,\quad d(x,y)\leq Cd(y,x)$$
In fact, we could prove Theorem \ref{th1} only supposing that $d$ is
weakly symmetric. But for simplicity, we only prove it in the true
metric setting.

Let us start with some generalities. First, note that up to
replacing $U$ by $U^m$, for some fixed $m>0$, we can assume that
$1\in U$, so that the sequence $(U^n)_{n\in\N}$ is nondecreasing.
Moreover, we have the following simple fact.
\begin{prop}\label{prop}
Let $G$ be a cglc group and let $U$ and $V$ be two compact sets such
that $U$ generates $G$ and contains $1$. Then there exists $m\in
\N^*$ such that, for all $n\geq m$, $V\subset U^{n}$.
\end{prop}

\bpr First, note that by a simple Baire argument, $U^n$ contains a
nonempty open set $\Omega$ for $n$ big enough. On the other hand,
for $n$ big enough, $U^n$ contains the inverse of a given element of
$\Omega$. Thus, $U^{n+1}$ contains an open neighborhood of $1$. Let
$\Omega.x_i$ be a finite covering of $V$. For $n$ big enough, we can
suppose that $x_i\in U^n$, so actually, $V\subset U^{2n+2}$. \epr

\begin{defn}
Let $G$ be a cglc group.
\begin{itemize}
\item We say that $G$ has polynomial growth if there exist a
compact generating set $U$, $D>0$ and a constant $C\geq 1$ such that
$$\mu(U^n)\leq Cn^D.$$
\item We say that $G$ has strictly polynomial growth if there exist a compact generating set $U$, a nonnegative number $d$, and a constant $C=C(U)\geq 1$ such that
\begin{equation}\label{eq}
C^{-1}n^d\leq \mu(U^n)\leq Cn^d.
\end{equation}
\end{itemize}
\end{defn}
Note that by Proposition~\ref{prop}, if $G$ has strict polynomial
growth, then the number $d$ does not depend on $U$, provided that
$\mu(U)\neq 0$. We call it the growth exponent of $G$.

\begin{thm}\cite{Gro,Gui,Los,Most,Wang}\label{thpoly}
Let $G$ be a cglc group of polynomial growth, then it has strictly
polynomial growth with integer exponent.
\end{thm}
Let us recall briefly how this result was proved. Using a structure
Theorem due to Wang \cite{Wang} and Mostow \cite{Most},
Guivarc'h~\cite[Corollary~III.3]{Gui} proved Theorem~\ref{thpoly}
for cglc solvable groups. Then a major step has been achieved by
Gromov~\cite{Gro}, who proved that a finitely generated group with
polynomial growth is virtually a lattice in some nilpotent connected
Lie group. Generalizing Gromov's approach, Losert~\cite{Los,Los'}
proved a similar statement for general cglc groups with polynomial
growth. According to Losert~\cite{Los}, $G$ is quasi-isometric to a
solvable cglc group $S$, and hence has strictly polynomial growth
with integer exponent.

\

In the group setting, we obtain a slightly improved version of
Theorem \ref{th1}.

\begin{thm}\label{th3}
Let $G$ be a cglc group with polynomial growth. Consider a sequence
$(U_n)_{n\in \N}$ of measurable subsets such that there exists two
generating compact sets $K, K'$ such that, for all $n\in \N,$
$$K\subset U_n\subset K'.$$
Write
$$N_n=U_0\ldots U_{n-1}U_n \quad \forall n\in \N$$
Then, there exist $\delta>0$ and a constant $C\geq 1$ such that
$$\mu\left(N_{n+1}\smallsetminus N_{n}\right)\leq Cn^{-\delta}\mu(N_n)\quad \forall n\in \N^*.$$
In particular, the sequence $(N_n)_{n\in\N}$ is F\o lner.
\end{thm}

The following corollary is a also a corollary of Theorem~\ref{th1}.

\begin{cor}\label{th2}
Let $G$ be a cglc group with polynomial growth, and $U$ be a compact
generating set of $G$. Then, there exist $\delta>0$ and a constant
$C\geq 1$ such that
$$\mu\left(U^{n+1}\smallsetminus U^n\right)\leq Cn^{-\delta}\mu(U^n)\quad \forall n\in \N^*.$$
In particular, the sequence $(U^n)_{n\in N}$ is F\o lner.
\end{cor}

In fact, we do not use the full contents of Theorem~\ref{thpoly}. We
only need the fact that $G$ satisfies a doubling property: there
exists of a constant $C=C(U)\geq 1$ such that
$$\mu(U^{2n})\leq C\mu(U^n) \quad \forall n\in \N.$$
It clearly results from Strictly Polynomial Growth. On the other
hand, Doubling Property implies trivially Polynomial Growth.
Unfortunately, there exist no elementary proofs of the converses,
which require to prove Theorem~\ref{thpoly}.

\section{Consequences in ergodic theory}
Let $G$ be a locally compact second countable (lcsc) group, $X$ a
standard Borel space on which $G$ acts measurably by Borel
automorphisms. Let $m$ be a $G$-invariant probability measure on $X$
($(X,m)$ is called a Borel probability $G$-space). The $G$-action on
$X$ gives rise to a strongly continuous representation $\pi$ of $G$
as a group of isometries of the Banach space $L^p(X)$ for $1\leq
p<\infty$, given by $\pi(g)f(x)=f(g^{-1}x)$. For any Borel
probability measure $\beta$ on $G$, and given some $p\geq 1$, we can
consider the averaging operator given by
$$\pi(\beta)f(x)=\int_G f(g^{-1}x)d\beta(g), \quad \forall f\in L^p(X).$$
Let $(\beta_n)$ be a sequence of probability measures on $G$. We say
that $(\beta_n)$ satisfies a pointwise ergodic theorem in $L^p(X)$
if
$$\lim_{n\rightarrow\infty}\pi(\beta_n)f(x)=\int_Xf dm$$
for almost every $x\in X$, and in the $L^p$-norm, for all $f\in
L^p(X)$, where $1\leq p<\infty$. Let $\mu$ be a Haar measure on $G$.
We will be interested in the case when $\beta$ is the normalized
average on a set of finite measure $N$ of $G$.
\begin{defn}[Regular sequences] A sequence of sets of finite
measure $N_k$ in $G$ is called regular if
$$\mu(N_k^{-1}.N_k)\leq C\mu(N_k).$$
\end{defn}

Let us recall the following general result (also proved in the
recent survey of Amos Nevo \cite{Amos}).

\begin{thm}\label{bigth}\cite{Bew,Chat,Emer,Temp}
Assume $G$ is an amenable lcsc group, and $(N_n)_{n\in \N}$ is an
increasing left F\o lner regular sequence, with $\cup_{n\in \N}
N_n=G$. Then, the sequence $(\beta_n)_{n\in \N}$ (associated to
$(N_k)$) satisfies the pointwise ergodic theorem in $L^p(X)$, for
every Borel probability $G$-space $(X,m)$ and every $1\leq
p<\infty$.
\end{thm}

Now, let us focus on the case when $G$ is a locally compact,
compactly generated group of polynomial growth. Consider a sequence
$(U_n)_{n\in \N}$ satisfying the hypothesis of Theorem \ref{th3}.
According to Theorem \ref{th3} and Proposition \ref{prop}, the
sequence $N_n=U_0.U_1\ldots U_n$ satisfies the hypothesis of Theorem
\ref{bigth}. So we get the following corollary.

\begin{thm}\label{thergo}
Let $G$ be a cglc group of polynomial growth. Consider a sequence
$(U_n)_{n\in \N}$ of measurable subsets such that there exist two
generating compact subsets $K, K'$ such that, for all $n\in \N$
$$K\subset U_n\subset K'.$$
Write $N_n=U_0U_1\ldots U_n$. Then, the sequence $(\beta_n)_{n\in
\N}$ (associated to $(N_n)_{n\in \N}$) satisfies the pointwise
ergodic theorem in $L^p(X)$, for every Borel probability $G$-space
$(X,m)$ and every $1\leq p<\infty$.
\end{thm}

\section{Remarks and questions}

In this section, we address a (non-extensive) list of remarks and
problems related to the subject of this paper.

\begin{que}
Is the Greenleaf localisation conjecture true for groups with
subexponential growth?
\end{que}

\begin{que}{\bf (Groups with exponential growth.)} Let $G$ be a finitely
generated group with exponential growth and let $U$ be a finite
generating subset. Does there exist a constant $c>0$ such
that\footnote{An erroneous proof of this fact is written in
\cite{P'}.}
$$\mu(U^{n+1}\smallsetminus U^n)\geq c\mu(U^n)?$$
\end{que}

\begin{que}{\bf (Asymptotic isoperimetry.)} Let $G$ be a locally compact,
compactly generated group and let $U$ be a compact generating
neighborhood of $1$. If $A$ is a subset of $G$, we call boundary of
$A$ and denote by $\partial A$ the subset $UA\smallsetminus A$. Let $\mu$
be a Haar measure on $G$. Recall the definition of the monotone
isoperimetric profile of $G$ (see \cite{PitSal,tess'})
$$I^{\uparrow}(t)=\inf_{\mu(A)\geq t}\mu(\partial A)/\mu(A)$$
where $A$ runs over measurable subsets of finite measure of $G$. We
can also define a (monotone) profile relatively to a family
$\mathbf{A}$ of subsets of $G$
$$I_{\mathbf{A}}^{\uparrow}(t)=\inf_{\mu(A)\geq t, A\in \mathbf{A}}\mu(\partial A)/\mu(A).$$
We say \cite{tess'} that the family $\mathbf{A}$ is asymptotically
isoperimetric if  $I_{\mathbf{A}}^{\uparrow}\preceq I^{\uparrow}$.

By a theorem of Varopoulos (\cite{Varo} \cite{Coul}), $G$ has
polynomial growth of exponent $d$ if and only if
$I^{\uparrow}(t)\approx t^{(d-1)/d}$.

An interesting question is for which groups do we have
$I_{(U^n)_{n\in \N}}^{\uparrow}\preceq I^{\uparrow}$?
\end{que}
It is true for groups of polynomial growth as shown by the following
proposition, valid for a general doubling metric measure space.
\begin{prop}\label{profil}
Let $X$ be a doubling metric measure space. There exists a sequence
$(r_i)_{i\in \N}$ such that $2^i\leq r_i\leq 2^{i+1}$ and such that
$$\forall i\in \N, \forall x\in X, \quad \mu(S(x,r_i))\leq C\mu(B(x,r_i)/r_i.$$
\end{prop}
In particular, if $G$ has polynomial growth of exponent $d$, and if
$U$ is a compact generating set of $G$, then there exists a
subsequence $n_i$ such that $2^i\leq n_{i+1}\leq 2^{i+1}$ and such
that
\begin{equation}\label{eq5}
\mu(U^{n_i+1}\smallsetminus U^{n_i})\leq Cn_i^{(d-1)/d}.
\end{equation}

\noindent{\bf Proof of Proposition~\ref{profil}.} First, remark that
$$\forall n<m\in \N, \quad S(x,n)\cap S(x,m)=\emptyset$$
and that
$$\cup_{k=1}^{2^i}S(x,2^i+k)\subset B(x,2^{i+1}),$$
so that
$$2^i\inf_{1\leq k\leq 2^i}\mu(S(x,2^i+k))\leq \mu(B(x,2^{i+1}))$$
and finally, one can conclude thanks to doubling property.\epr

\begin{rem}
In a general graph with strict polynomial growth of exponent $d$,
the profile may sometimes be much smaller than $t^{(d-1)/d}$, so
Proposition~\ref{profil} does not necessarily imply that balls are
asymptotically isoperimetric, i.e.
$I_{(B(x,k))_{x,r}}^{\uparrow}\preceq I^{\uparrow}.$ This issue is
studied quite extensively in \cite{tess'}.
\end{rem}

\begin{que}
One can make the last question more precise by asking if
$I_{(U^k)}^{\uparrow}\preceq I^{\uparrow}$ implies that $G$ has
polynomial growth? Subexponential growth?
\end{que}

\begin{que}\label{ques}
A very natural question is: does (\ref{eq5}) hold for any integer
$n$ (when $G$ has polynomial growth of exponent $d$), or
equivalently, is there a constant $C < \infty$ such that:
\begin{equation}\label{isop}
\forall n\in \N,\quad \mu(U^{n+1}\smallsetminus U^{n})\leq
C\frac{\mu(U^n)}{n}?
\end{equation}
\end{que}
This question is motivated by the following observation.
\begin{prop}\label{abelian}
Let $G$ be a cglc abelian group and let $U$ be a compact generating
set of $G$. Then, (\ref{isop}) holds.
\end{prop}
\noindent{\bf Sketch of the proof.} First, note that it is an easy
fact when $G=\R^d$ (the adaptation to $\Z^d$ is left to the reader):
if $K$ is convex, it is trivial (since $K+K=2.K$); then show that
$\hat{K}^n\subset K^{n+k}$ where $\hat{K}$ denotes the convex hull
of $K$, and where $k$ is a positive integer smaller than $d+1$ times
the diameter of $K$. On the other hand, a cglc Abelian group $G$ is
isomorphic to a direct product $K\times \R^a\times \Z^b$, with
$a,b\in \N$, and $K$ being a compact group. \epr

\begin{rem}
Question~\ref{ques} is also natural in the context of doubling
graphs. In this setting, the question becomes: does there exist a
constant $C$ such that for all $x\in X$, and all $n\geq 1$,
$$\frac{|S(x,n)|}{|B(x,n)|}\leq \frac{C}{n}?$$
But as mentioned in the introduction, the answer is no in a very
strong sense since in \cite[Theorem~4.9]{tess'}, we construct a
graph quasi-isometric to $\Z^2$ that does not satisfy this property.
\end{rem}

\section{Proofs}

We will start proving Theorem \ref{th1} which is our ``more general
result". Nevertheless, Theorem \ref{th3} is not an immediate
consequence of the group version of Theorem \ref{th1}, that is,
Corollary \ref{th2}. So for the convenience of the reader, we will
give a proof of Corollary \ref{th2} using notation adapted to the
group setting, and then give the additional argument which is needed
to obtain Theorem \ref{th3}.

\subsection{A preliminary observation}

\

The following observation is one of the main ingredients of the
proofs.
\begin{lem}\label{th}
Let $X=(X,\mu)$ be a measured space. Let us consider an increasing
sequence $(A_n)_{n\in \N}$ of measurable subsets of $X$. Define
$C_{n,n+k}=A_{n+k}\smallsetminus A_n$. We suppose that $\mu(A_n)$
is finite and unbounded with respect to $n\in \N$. Let us suppose
that there exists a constant $\alpha>0$ such that, for all
integers $k\leq n$,
\begin{equation}\label{equ}
\mu(C_{n-k,n})\geq \alpha.\mu(C_{n,n+k}).
\end{equation}
Then, there exist $\delta>0$ and a constant $C\geq 1$ such that
$\forall n\geq 1$
$$\frac{\mu(C_{n-1,n})}{\mu(A_n)}\leq Cn^{-\delta}.$$
\end{lem}

\bpr Write $i_n=[\log_2n]$. For $i\leq i_n$, define
$b_{i}=\mu(C_{n-2^i,n})$. Note that
$$C_{n-2^{i},n}=C_{n-2^{i},n-2^{i-1}}\cup\ldots \cup C_{n-1,n}
\quad \forall i\leq i_n$$ and that the union is disjoint. So we have
$$b_i=\mu(C_{n-2^{i},n-2^{i-1}})+\ldots+\mu(C_{n-1,n}).$$ On the
other hand, by \ref{equ}
\begin{eqnarray*}
\mu(C_{n-2^{i},n-2^{i-1}})&  =   &
\mu(C_{n-2^{i-1}-2^{i-1},n-2^{i-1}})\\
                          & \geq &
\alpha.\mu(C_{n-2^{i-1},n-2^{i-1}+2^{i-1}})\\
                          &  =   &
\alpha.b_{i-1}.
\end{eqnarray*}
But note that
$$b_{i}=b_{i-1}+\mu(C_{n-2^{i},n-2^{i-1}})$$
So
$$b_{i}\geq (1+\alpha)b_{i-1}.$$
Therefore
$$b_{i}\geq (1+\alpha)^i\mu(C_{n-1,n}).$$
Thus, it comes
\begin{eqnarray*}
\mu(A_n)  & \geq & b_{i_n}\\
          & \geq &  (1+\alpha)^{i_n}\mu(C_{n-1,n})\\
          & \geq & (1+\alpha)^{\log_2n-1}\mu(C_{n-1,n})\\
          & \geq & \frac{1}{1+\alpha}n^{\log_2(1+\alpha)}\mu(C_{n-1,n}).
\end{eqnarray*}
So we are done. \epr

\subsection{The case of metric measured spaces: proof of the theorem \ref{th1}}\label{sec}

For all $x\in X$ and $r'>r>0$, write
$$C_{r,r'}(x)=B(x,r')\smallsetminus B(x,r),$$
and
$$c_{r,r'}(x)=\mu(C_{r,r'}(x)).$$
Thanks to Lemma~\ref{th}, we only need to prove that shells are
doubling, i.e. that there exists a constant $\alpha>0$ such that,
for all $x\in X$, and for any integers $n>k>10C$ (where $C$ is the
constant that appears in the definition of Property (M))
$$c_{n-k,n}(x)\geq \alpha c_{n,n+k}(x).$$
So it is enough to prove the following lemma:

\begin{lem}\label{lem2}
Let $(X,d,\mu)$ a doubling (M) space. Then, there exists $\alpha>0$
such that for all $x\in X$ and for all couples of integers $4C<k\leq
n$,
$$c_{n-k,n}(x)\geq \alpha.c_{n,n+k}(x).$$
\end{lem}
\bpr Let $y$ be in $C_{n,n+k}(x)$. Consider a finite chain $x_0=y,
x_1,\ldots,x_m=x$ such that for $0\leq i<m$
$$d(x_i,x_{i+1})\leq C;$$
and
$$d(x_{i+1},x)\leq d(x_{i},x)-1.$$
Let $k_0$ be the smallest integer such that $x_{k_0}\in B(x,n-k/2)$.
Since $y\in C_{n,n+k}(x)$, $k_0$ exists and is less than $2k$.
Moreover, minimality of $k_0$ implies that $x_{k_0}\in
C_{n-k/2-C,n-k/2}(x)$. So we have
\begin{equation}\label{equ'}
d(y,C_{n-k/2-C,n-k/2}(x))\leq 2Ck.
\end{equation}

Let $(z_i)_i$ be a maximal family of $k$-separated
points\footnote{As the space is doubling, such a family exists and
is finite.} in $C_{n-k/2-C,n-k/2}(x)$. Provided for instance that
$k\geq 2C$, $C_{n-k/2-C,n-k/2}(x)$ is covered by the balls
$B(z_i,2k)$. Consequently, (\ref{equ'}) implies that the balls
$B(z_i,(2+2C)k)$ cover $C_{n,n+k}(x)$. On the other hand, if $k\geq
4C$, then for every $i$, $z_i$ belong to $C_{n-3k/4,n-k/2}$ and
hence the ball $B(z_i,k/4)$ is included in $C_{n-k,n}(x)$. Moreover,
these balls are disjoint. So we conclude by doubling property. \epr

\

\begin{rem}
Note that a lower bound on $k$ depending on $C$ is necessary because
otherwise, $C_{n-k,n}(x)$ could be empty for $k=1$. For instance
consider $\Z$ equipped with usual distance multiplied by $2$: it
satisfies Property (M) with $C=2$. In this case, $C_{n-1,n}(0)$ is
empty for odd $n$ although $C_{n,n+1}$ is not empty.
\end{rem}

\subsection{The case of groups: proof of Theorem \ref{th3}}

Note that to prove Theorem~\ref{th3}, we can assume that $U_n$
contains $1$, at least for $n$ large enough, which ensures that
$(N_n)$ is nondecreasing. Indeed, choose an integer $m$ such that
$K'\cup\{1\}\subset K^m$. Then, write $n=qm+r$ with $r<m$ and for
all $j\geq 1$, define
$$\tilde{U}_j=U_{(j-1)m+r+1}\ldots U_{jm+r}.$$
Define also $\tilde{U}_0=U_0\ldots U_r$. We therefore have
$$\tilde{N}_q=N_n=\tilde{U}_0\ldots \tilde{U}_q.$$
Finally, as $U_{n+1}\subset K'\subset \tilde{U}_{n+1}$, it suffices
to prove Theorem~\ref{th3} for the sequence $(\tilde{N_q})$.

Actually, instead of directly proving Theorem~\ref{th3}, we will
prove Corollary~\ref{th2} and then explain how the proof can be
generalized.

\

\noindent{\bf Proof of Corollary~\ref{th2}.} Let $G$ be a cglc group
of polynomial growth endowed with a Haar\footnote{Note that since
$G$ has subexponential growth, it is unimodular, so that the Haar
measure is left and right invariant.} measure $\mu$. Let $U$ be a
compact generating subset containing $1$. Recall that this implies
that the sequence $U^n$ is nondecreasing. Let us write
$$C_{n,n+k}=U^{n+k}\smallsetminus U^n, \quad \forall n,k\in \N$$
and
$$c_{n,n+k}=\mu(C_{n,n+k}).$$
Recall that we want to find a constant $\alpha$ such that
$c_{n-k,n}\geq \alpha.c_{n,n+k}$ for $k$ large enough. To simplify
notation, let us assume that $k$ is a positive multiple\footnote{If
$k$ is not a multiple of $k$, one has to assume at least that $k\geq
4$ and to replace everywhere in the proof $k/4$ and $k/2$ by their
integer parts.} of $4$.

First, we have
\begin{clai}\label{eq1}
$$C_{n,n+k}\subset C_{n-k/2,n-k/2+1}U^{2k}.$$
\end{clai}
\bpr Indeed, let $y$ be in $C_{n,n+k}$, and let $(y_1,\ldots
,y_{n+j})$ be a minimal sequence of elements of $U$ such that
$y=y_{1}\ldots y_{n+j}$. By definition of $C_{n,n+k}$ and by
minimality, we have $1\leq j\leq k$. Moreover, it is easy to see
that minimality also implies
$$y_{1}\ldots y_{n-k/2+1}\in
C_{n-k/2,n-k/2+1}.$$ So
$$y\in C_{n-k/2,n-k/2+1}y_{n-k/2+2}\ldots
y_{n+j}\subset C_{n-k/2,n-k/2+1}U^{2k}$$ and we are done.\epr

\

On the other hand, we have
\begin{clai}
$$C_{n-k/2,n-k/2+1}U^{k/4}\subset C_{n-k,n}.$$
\end{clai}
\bpr Since $U$ contains $1$, we have
$$C_{n-k/2,n-k/2+1}U^{k/4}\subset U^{n-k/2+k/4+1}\subset U^n.$$
Besides, let $x\in C_{n-k/2,n-k/2+1}U^{k/4}$, so that
$$x=yu_1\ldots u_{k/4}.$$
If we had $x\in U^{n-k},$ then, it would imply that $y\in
U^{n-k+k/4}\subset U^{n-k/2}$: absurd. So $x\in C_{n-k,n}.$\epr

\

Now, let $(x_i)$ be a maximal family of points of
$C_{n-k/2,n-k/2+1}$ such that $x_iU^{k/4}\cap x_jU^{k/4}=\emptyset$
for $i\neq j$. By maximality of $(x_i)$, we have
$$C_{n-k/2,n-k/2+1}\subset\cup_i
x_iU^{k/4}U^{-k/4}.$$ So by Claim~\ref{eq1}, we get
\begin{equation}\label{eq2}
C_{n,n+k}\subset \cup_i x_iU^{k/4}U^{-k/4}U^{2k}
\end{equation}
Let $S$ be a symmetric compact neighborhood of $1$ containing $U^3$.
Then, since $U^{k/4}U^{-k/4}U^{2k}$ is included in $S^k$,
Theorem~\ref{thpoly} implies that there exists a constant $C<\infty$
such that for $k$ large enough,
\begin{equation}\label{eq3}
\mu\left(x_iU^{k/4}U^{-k/4}U^{2k}\right)\leq
C\mu\left(x_iU^{k/4}\right).
\end{equation}
Thus, since the $x_iU^{k/4}$ are disjoint and included in
$C_{n-k,n}$, we get
\begin{equation}\label{eq4}
c_{n-k,n}\geq \sum_i \mu\left(x_iU^{k/4}\right).
\end{equation}
Finally, using (\ref{eq2}), (\ref{eq3}) and (\ref{eq4}), we deduce
$$c_{n-k,n}\geq C^{-1}c_{n,n+k}.\;\blacksquare$$

\

\noindent{\bf Proof of Theorem \ref{th3}} The only significant
modification we have to do in order to prove Theorem \ref{th2}
concerns Claim~\ref{eq1}. Actually, we have to show a kind of
Property (M) adapted to this context.

\begin{clai}\label{eq6}
There exists $j_0\in \N$, such that for every $n\in \N$ and every
$x\in U_0\ldots U_{n+k} \smallsetminus U_0\ldots  U_n$, we have
$$x\in K'^k\left(U_0\ldots U_n\smallsetminus U_{0}\ldots U_{n-kj_0}\right).$$
\end{clai}
\bpr Since $K'$ contains $U_{i}$ for every $i\in \N$, we have
$$x\in U_0 \ldots U_nK'^k.$$
On the other hand, let $q$ be an integer such that
$$x\in U_{0}\ldots U_{n-q}K'^k.$$
Then, let $j_0$ be such that $K'\subset K^{j_0}$ (see
Proposition~\ref{prop}). Since $K\subset U_i$ for every $i$, it
comes
$$x\in U_0\ldots U_{n-q+kj_0}.$$
But this implies $q<kj_0$, so we are done. \epr

\

Let us finish the proof of Theorem~\ref{th3}. Write
$$C_{n,n+k}=U_0\ldots U_{n+k}\smallsetminus U_0\ldots U_n.$$
According to Claim~\ref{eq6}, we have
$$C_{n,n+k}\subset C_{n-kj_0,n} K'^k$$
for every $k<n/j_0$.

Using the same arguments as in the proof of Theorem~\ref{th2}, we
get $$c_{n,n-j_0k}\geq \alpha.c_{n,n+k},$$ and we conclude thanks to
Lemma~\ref{th}. \epr

\bigskip
\begin{center}
{\bf Acknowledgements}
\end{center}

I am grateful to A. Nevo for having pointed to me the ergodic
consequences of the present result. I also thank P. de la Harpe, T.
Coulhon, Y. de Cornulier, E. Breuillard and the referee for their
valuable comments and corrections.

\bibliographystyle{amsplain}

\bigskip
\footnotesize
\noindent Romain Tessera\\
Université de Cergy-Pontoise\\
E-mail: \url{tessera@clipper.ens.fr}\\

\end{document}